\documentclass{aip-cp}

\usepackage[numbers]{natbib}
\usepackage{rotating}
\usepackage{graphicx}

\begin{document}

\newtheorem{theorem}{Theorem}
\newtheorem{definition}{Definition}
\newtheorem{corollary}[theorem]{Corollary}
\newtheorem{lemma}[theorem]{Lemma}
\newtheorem{proposition}[theorem]{Proposition}

% Title portion
\title{Relativistic formulation of abstract evolution equations}
\title{Relativistic formulation of abstract evolution equations}

\author[aff1,aff2]{Yoritaka Iwata\corref{cor1}}
%\eaddress[url]{http://www.aip.org}
%\eaddress{shirakawa-kiyoaki@kakaku.com}
%%%

\affil[aff1]{Institute of Innovative Research, Tokyo Institute of Technology, Japan}
\affil[aff2]{Department of Mathematics, Shibaura Institute of Technology, Japan}
%\affil[aff3]{You would list an author's second affiliation here.}
\corresp[cor1]{E-mail: iwata\_phys@08.alumni.u-tokyo.ac.jp}
%\authornote[note1]{The autuhor is grateful to Prof. Emeritus Hiroki Tanabe of Osaka University..}

\maketitle

\begin{abstract}
The relativistic formulation of abstract evolution equations is introduced. The corresponding logarithmic representation is shown to exist without assuming the invertible property of evolution operators. Consequently, by means of the logarithmic representation of operators, nonlinear and discrete properties are shown to be valid to the infinitesimal generators of relativistic abstract evolution equations. 
\end{abstract}

\section{INTRODUCTION}
The relativistic formulation of abstract evolution equations is introduced by Ref.~\cite{18iwata} in order to establish an abstract version of the Cole-Hopf transform in Banach spaces and to explain the nonlinear relation between the evolution operator and its infinitesimal generator. 
In the present article the logarithmic representation of infinitesimal generator is obtained for the relativistic form of abstract evolution equations.
It is a generalization of the logarithmic representation for the infinitesimal generators of invertible evolution families.

Let $X$ and $B(X)$ be a Banach space and a space of bounded linear operators on $X$, respectively.
The norm notations $\| \cdot \|_X$ and $\| \cdot \|_{B(X)}$ are used for the norms equipped with $X$ and $B(X)$, respectively.
For a positive and finite $T$, let $t$ and $s$ satisfy $-T \le t,s \le T$, $Y$ be a dense Banach subspace of $X$, and the topology of $Y$ be stronger than that of $X$.
The space $Y$ is assumed to be $U(t,s)$-invariant.
The elements of evolution family $\{U(t,s) \}_{-T \le t, s  \le T}$ are assumed to be mappings: $(t,s) \to U(t,s)$ satisfying the strong continuity for $-T \le t, s  \le T$.
The semigroup properties:
\begin{equation} \label{sg1} 
U(t,r)~ U(r,s) = U(t,s), 
\end{equation}
and
\begin{equation}  \label{sg2} 
U(s,s) = I, 
\end{equation}
are assumed to be satisfied, where $I$ denotes the identity operator of $X$.
Although both $U(t,s)$ and $U(s,t)$ have been assumed to be well-defined to satisfy $U(s,t) ~ U(t,s) = U(s,s) = I$
in the preceding articles~\cite{18iwata}, the invertible property of $U(t,s)$ is not assumed in this article.
Consequently, $U(t,s)$ is the two-parameter semigroup (for a textbook, see \cite{79tanabe}) on a finite or infinite dimensional Banach space $X$.
Note that $U(t,s)$ can be either linear or nonlinear semigroup.

\section{PRECEDING RESULTS}
\begin{figure}[t] 
\includegraphics[width=60mm]{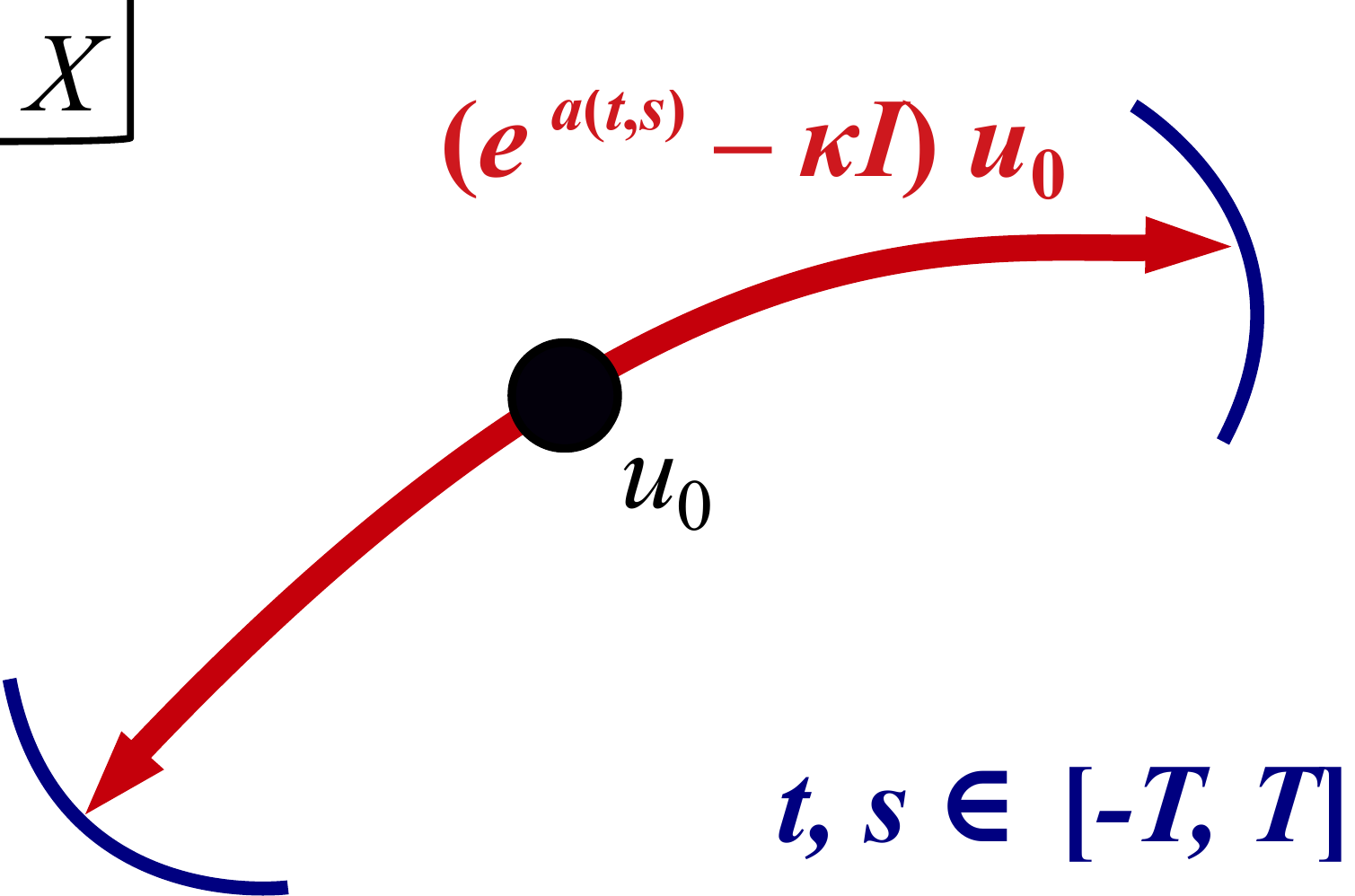} \vspace{6mm}
\caption{ \label{fig1} 
(Color online) The regularized trajectory $(e^{a(t,s)}-\kappa I) u_0 \in X$ generated by the alternative infinitesimal generator $\partial_t a(t,s)$, where $u_0 \in X$ denotes a certain element in $X$.
The regularized trajectory is a local representation which is valid only within a finite interval of $t$.
 }
\end{figure}

The pre-infinitesimal generator of $U(t,s)$ is introduced based on Ref.~\cite{17iwata-1}.
The pre-infinitesimal generator is a generalized concept of infinitesimal generator, but it is not necessarily an infinitesimal generator without assuming a dense property of domain space $Y$ in $X$.
That is, only the exponentiability with a certain ideal domain is valid to the pre-infinitesimal generators. 
For $u_s \in Y$ and $-T \le t, s  \le T$, the linear operator $A(t): Y  ~\to~  X$ is called the pre-infinitesimal generator if a weak limit
\begin{equation} \label{pe-group}
A(t) u_s := \mathop{\rm wlim}\limits_{h \to 0}  h^{-1} (U(t+h,t) - I) u_s
\end{equation}
exists.
Using the pre-infinitesimal generator, the definition of $t$-differential of $U(t,s)$ in a weak sense
\begin{equation} \label{de-group} \begin{array}{ll} 
\partial_t U(t,s)~u_s = A(t) U(t,s) ~u_s
\end{array}  \end{equation}
follows.
Accordingly the operator $A(t)$ is generally unbounded in $X$, and Eq.~(\ref{de-group}) is regarded as the abstract evolution equation evolving for $t$-direction (for example, see \cite{79tanabe}).
%%%%%% 
Consequently, for $-T \le t \le T$, $A(t)$ is represented by means of the logarithm function~\cite{17iwata-1}; there exists a certain complex number $\kappa \ne 0$ such that
\begin{equation}  \begin{array}{ll}
\label{logex}  A(t) ~ u_s =  (I+ \kappa U(s,t))~ \partial_{t} a(t,s)  u_s, \vspace{2.5mm} \\
 \partial_{t} a(t,s) = \partial_{t} {\rm Log}  (U(t,s) + \kappa I),
\end{array} \end{equation}
where ${\rm Log}$ denotes the principal branch of logarithm, and $\partial_{t} a(t,s) $ is the alternative infinitesimal generator~\cite{17iwata-3} to $A(t)$. 
By integrating the second equation of Eq.~(\ref{logex}), $ \exp [a(t,s)] =  U(t,s) + \kappa I $ is obtained.
It shows a correspondence between $\exp [a(t,s)]$ and $U(t,s)$.
There are some differences between $a(t,s) $ and $A(t)$. 
For example, $a(t,s) $ is always a bounded operator, while $A(t)$ is not necessarily a bounded operator. 
This is an advantage of introducing the alternative infinitesimal generator providing a regularization.
It leads to the local-in-$t$ regularized trajectory (Fig.~\ref{fig1}).

\section{MAIN RESULT}
Let the standard space-time variables $(t,x,y,z)$ be denoted by $(x^0,x^1,x^2,x^3)$ receptively.
It is further possible to generalize space-time variables to $(x^0,x^1,x^2,x^3, \cdots, x^{n})$ being valid to general $(n+1)$-dimensional space-time. 
In spite of the standard treatment of abstract evolution equations, the direction of evolution does not necessarily mean time-variable $t=x^0$ in the relativistic treatment of the evolution equations.
The condition to obtain the logarithmic representation is stated as follows.

\begin{definition}
For an evolution family of operators $\{U(x^i,\xi^i)\}_{-L \le x^i,\xi^i \le L}$ in a Banach space $X_i$, let  $K(x^i):Y_i \to X_i$ be the infinitesimal generator of $U(x^i,\xi^i)$, where $Y_i$ is a dense subspace of $X_i$.
The relativistic form of abstract evolution equations is defined as
\begin{equation} \label{k-eq1}  \begin{array}{ll} 
\partial_{x^i} U(x^i, \xi^i)~u({\xi^i})  = K(x^i) U(x^i,\xi^i) ~u({\xi^i}) , \vspace{2.5mm} \\
 u({\xi^i}) = u_{\xi^i}
 \end{array}  \end{equation}
in $X_i$, where $X_i$ is a functional space consisting of functions with variables $x^j$ with $0 \le j \le n$ skipping only $j=i$.
Consequently, the unknown function is represented by $u(x^i)  =  U(x^i,\xi^i) ~u_{\xi^i}$ for a given initial value $u_{\xi^i} \in X_i$.
\end{definition}

Let $\partial_{x^i} U(x^i, \xi^i)~u_{\xi^i} = K(x^i) U(x^i,\xi^i) ~u_{\xi^i}$ evolving for $i$-direction be represented by $\partial_{x^k} V(x^k, \xi^k)~v_{\xi^k} =  {\mathcal K}(x^k) V(x^k,\xi^k) ~v_{\xi^k}$ in a certain direction $k$.
For $k \ne i$, we consider
\begin{equation} \label{k-eq2}  \begin{array}{ll} 
\partial_{x^k} V(x^k, \xi^k)~v({\xi^k}) = {\mathcal K}(x^k) V(x^k,\xi^k) ~v({\xi^k}), \vspace{2.5mm} \\
 v({\xi^k}) = v_{\xi^k}
 \end{array}  \end{equation}
in $X_k$.
It is remarkable that even if the evolution operator $U(x^i, \xi^i)$ and its infinitesimal generator exist, $V(x^k, \xi^k)$ and its infinitesimal generator do not necessarily exist.
Those existence should be individually examined for each direction (for an example, see Ref~\cite{18iwata}).
If $U(x^i, \xi^i)$, $V(x^k, \xi^k)$ and those infinitesimal generators exist, $u(x^i)$ in Eq.~(\ref{k-eq1}) and $v(x^k)$ in Eq.~(\ref{k-eq2}) satisfy the same evolution equation, where the detailed conditions such initial and boundary conditions can be different depending on the settings of $X_i$ and $X_k$.  
The propagation of singularity should be different if the evolution direction is different. 
For Eqs.~(\ref{k-eq1}) and (\ref{k-eq2}), the evolution direction is not limited to $x^0$.
This gives a reason why the formulation shown in Eq.~(\ref{k-eq1}) is called the relativistic form of abstract evolution equations. 

\begin{theorem}
Let $i$ denote any direction satisfying $0 \le i \le n$.
Let $x^i$ and $\xi^i$ satisfy $-L \le x^i,\xi^i \le L$, and $Y_i$ be a dense subspace of a Banach space $X_i$.
A two-parameter evolution family of operators $\{U(x^i,\xi^i)\}_{-L \le x^i,\xi^i \le L}$ satisfying Eqs.~(\ref{sg1}) and (\ref{sg2}) is assumed to exist in a Banach space $X_i$.
Under the existence of the infinitesimal generator $K(x^i):Y_i \to X_i$ of $U(x^i,\xi^i)$ for the $x^i$ direction, let $U(x^i,\xi^i)$ and $K(x^i)$ commute.
The logarithmic representation of infinitesimal generator is obtained; there exists a certain complex number $\kappa \ne 0$ such that
\begin{equation} \label{logex2} \begin{array}{ll}
 K(x^i) ~ u_{\xi^i} =
 \{  I - \kappa (U(x^i,\xi^i) + \kappa I)^{-1} \}^{-1} ~ \partial_{x^i} {\rm Log} ~ (U(x^i,\xi^i) + \kappa I) ~ u_{\xi^i},
\end{array} \end{equation}
where $u_{\xi}$ is an element in $Y_i$, and $\kappa$ is taken from the resolvent set of $U(x^i,\xi^i)$.
Note that $U(x^i,\xi^i)$ is not assumed to be invertible.
\end{theorem}
%%%%%%
{\bf Proof} \quad
Different from the proof in the previous research \cite{17iwata-1}, here the similar statement is proved without assuming the invertible property of $U(x^i,\xi^i)$
The key point is that $ (U(x^i,\xi^i) + \kappa I)^{-1} $ exists for a certain $\kappa \in {\mathbb C}$, even if $ U(x^i,\xi^i)^{-1} $ does not exist.
In particular, the obtained representation is more general compared to the one obtained in Ref.~\cite{17iwata-1}.
For any $U(x^i,\xi^i)$, operators $ {\rm Log} ~ (U(x^i,\xi^i) + \kappa I)$ and $ {\rm Log} ~ (U(x^i+h,\xi^i) + \kappa I)$ are well defined for a certain $\kappa$.
%%%
The $x^i$-differential in a weak sense is formally written by
\begin{equation} \label{difference0} \begin{array} {ll} 
\mathop{\rm wlim}\limits_{h \to 0}  \frac{1}{h} \{ {\rm Log} ~(U(x^i+h,\xi^i)+\kappa I) - {\rm Log} ~(U(x^i,\xi^i)+ \kappa I) \}   \vspace{1.5mm} \\
 = \mathop{\rm wlim}\limits_{h \to 0}  \frac{1}{2 \pi i}
 \int_{\Gamma} {\rm Log} \lambda   %\vspace{1.5mm} \\ 
 ~ \{ (\lambda I - U(x^i+h,\xi^i)-\kappa I )^{-1} \frac{U(x^i+h,\xi^i)-U(x^i,\xi^i)}{h} (\lambda I - U(x^i,\xi^i) - \kappa I )^{-1} \} d \lambda   
\end{array} \end{equation} 
where $\Gamma$, which is taken independent of $x^i$, $\xi^i$ and $h$ for a sufficiently large certain $\kappa$, denotes a circle in the resolvent set of both $U(t,s)+ \kappa I$ and $U(t+h,s)+\kappa I$.
There are three steps to establish the validity of logarithmic representation in the relativistic form.

In the first step, the continuity of the mapping $x^i \to (\lambda I -  U(x^i,\xi^i) - \kappa I)^{-1}$ as for the strong topology is shown.
A part of the integrand of Eq.~(\ref{difference0}) is estimated as
\begin{equation} \label{intee} \begin{array} {ll}
\quad \| \{ (\lambda I -  U(x^i+h,\xi^i) -\kappa I )^{-1} \frac{ U(x^i+h,\xi^i)- U(x^i,\xi^i)}{h} (\lambda I -  U(x^i,\xi^i)- \kappa I )^{-1} \} v \|_{X_i} \vspace{1.5mm} \\
  \le  \| \{ (\lambda I -  U(x^i+h,\xi^i) - \kappa I )^{-1} \|_{B(X)}  %\vspace{1.5mm} \\
 \|  \frac{U(x^i+h,\xi^i)- U(x^i,\xi^i)}{h} (\lambda I -  U(x^i,\xi^i) - \kappa I )^{-1} \} v \|_{X_i}, 
\end{array} \end{equation}
for $v \in X_i$. 
The operator boundedness 
\[ \begin{array} {ll}
 \| \{ (\lambda I -  U(x^i+h,\xi^i) - \kappa I )^{-1} \|_{B(X_i)}  < \infty
\end{array} \]
is true, since $\lambda$ is taken from the resolvent set of $ U(x^i+h,\xi^i) - \kappa I$.
In the same way the operator $(\lambda I - U(x^i,\xi^i) - \kappa I )^{-1}$ is bounded on both $X_i$ and $Y_i$, where $U(x^i,\xi^i)$ is strongly continuous mapping with respect to the variable $x^i$ on both $X_i$ and $Y_i$.
Then the continuity of the mapping $x^i \to (\lambda -  U(x^i,\xi^i) - \kappa )^{-1}$ as for the strong topology follows:
\[ \begin{array} {ll}
 \|  (\lambda I -  U(x^i+h,\xi^i) - \kappa I )^{-1} - (\lambda I  -  U(x^i,\xi^i) - \kappa I)^{-1}  \|_{B(X_i)}   \vspace{1.5mm} \\
 \le  \| (\lambda I -  U(x^i+h,\xi^i) - \kappa I)^{-1} \|_{B(X_i)}   %\vspace{1.5mm} \\
 \|(  U(x^i+h,\xi^i) - U(x^i,\xi^i)) (\lambda I - U(x^i,\xi^i)-\kappa I )^{-1}  \|_{B(X_i)}.
\end{array} \]

In the second step, the uniform convergence is shown.
The latter part of the right hand side of Eq.~(\ref{intee}) is estimated as
\begin{equation} \label{unibound}  \begin{array}{ll}
 \left\| \frac{ U(x^i+h,\xi^i)- U(x^i,\xi^i)}{h} (\lambda I  -  U(x^i,\xi^i) -\kappa I )^{-1}  u  \right\|_{X_i} %
 \le  \frac{1}{|h|}  \int_{x^i}^{x^i+h} \|  K(\eta^i)   U(\eta^i,\xi^i) \|_{B(Y_i,X_i)}%\vspace{1.5mm}\\ 
   \| (\lambda I - U(x^i,\xi^i)- \kappa I )^{-1}\|_{B(Y_i)} \|u \|_{Y_i} ~ d\eta 
\end{array} \end{equation}
for $u \in Y_i$.
Because $ \|  K(\eta^i)   U(\eta^i,\xi^i) \|_{B(Y_i,X_i)} < \infty $ is true by assumption, the right hand side of Eq.~(\ref{unibound}) is finite.  
Equation~(\ref{unibound}) shows the uniform boundedness with respect to $h$, then the uniform convergence ($h \to 0$) of Eq.~(\ref{difference0}) follows.
Consequently, the weak limit $h \to 0$ for the integrand of Eq. (\ref{difference0}) is justified, as well as the commutation between the limit and the integral.

In the final step, the logarithmic representation is obtained.
According to Eq.~(\ref{difference0}), the interchange of the limit with the integral leads to
\[ \begin{array} {ll} 
\partial_t {\rm Log} (U(x^i,\xi^i) + \kappa I) ~ u 
=  \frac{1}{2 \pi i} \int_{\Gamma} d\lambda  %\vspace{1.5mm} \\    
 ~ \left[ ( {\rm Log} \lambda )  (\lambda I - U(x^i,\xi^i) - \kappa I )^{-1} 
 ~\mathop{\rm wlim}\limits_{h \to 0} \left ( \frac{U(x^i+h,\xi^i)-U(x^i,\xi^i)}{h} \right)
~ (\lambda I  - U(x^i,\xi^i)  - \kappa I )^{-1}    \right] ~u \\
\end{array} \] 
for $u \in Y_i$.
Because it is allowed to interchange $K(x^i)$ with $U(x^i,\xi^i)$,
\[  \begin{array}{ll}
\partial_t {\rm Log} (U(x^i,\xi^i) + \kappa I) ~ u %\vspace{1.5mm}\\
% = \frac{1}{2 \pi i} \int_{\Gamma} ({\rm Log} \lambda)
% (\lambda-U(t,s)-\kappa )^{-1} 
%  A(t) ~ U(t,s) ~ (\lambda -U(t,s)- \kappa )^{-1}  d \lambda ~ u \vspace{1.5mm}\\
 = \frac{1}{2 \pi i} \int_{\Gamma} ({\rm Log} \lambda) ~(\lambda I -U(x^i,\xi^i) - \kappa I)^{-2} ~ U(t,s) ~ d \lambda ~ K(x^i) ~ u
\end{array} \]
for $u \in Y_i$, where $\mathop{\rm wlim}\limits_{h \to 0} \left ({U(x^i+h,\xi^i)-U(x^i,\xi^i)} \right)/h$ means the pre-infinitesimal generator $K(x^i)$ itself.
A part of the right hand side is calculated as
\[ \begin{array}{ll}
 \quad  \frac{1}{2 \pi i} \int_{\Gamma}~ ({\rm Log} \lambda) ~(\lambda I -U(x^i,\xi^i) - \kappa I)^{-2} U(x^i,\xi^i) ~ d \lambda  %\vspace{1.5mm} \\
  = I- \kappa (U(t,s)+ \kappa I)^{-1},  
\end{array} \] 
due to the integration by parts, where the details of procedure is essentially the same as Ref.~\cite{17iwata-1}.
It leads to
\[ \begin{array}{ll}
 K(x^i) ~ u_{\xi^i} =
 \{  I - \kappa (U(x^i,\xi^i) + \kappa I)^{-1} \}^{-1} ~ \partial_{x^i} {\rm Log} ~ (U(x^i,\xi^i) + \kappa I) ~ u_{\xi^i},
%A(t) ~ u = \{I- \kappa (U(t,s)+\kappa I)^{-1}\}^{-1} ~ \partial_{t} {\rm Log} ~ (U(t,s) + \kappa I) ~ u  %\vspace{1.5mm} \\
%\quad = (U(t,s)+\kappa I) U(t,s)^{-1} ~ \partial_{t} {\rm Log} ~ (U(t,s) + \kappa I) ~ u  \vspace{1.5mm} \\
%\quad =  (I+\kappa U(s,t))~ \partial_{t} {\rm Log} ~ (U(t,s) + \kappa I) ~ u 
\end{array} \]
for $u_{\xi^i} \in Y_i$.
It is notable that $ (U(x^i,\xi^i) + \kappa I)^{-1}$ is always well defined for any $\kappa$ taken from the resolvent set of $U(x^i,\xi^i)$, even if $ U(\xi^i,x^i) = U(x^i,\xi^i)^{-1}$ does not exist.
\quad $\square$ \\

\begin{figure}[t] 
\includegraphics[width=48mm]{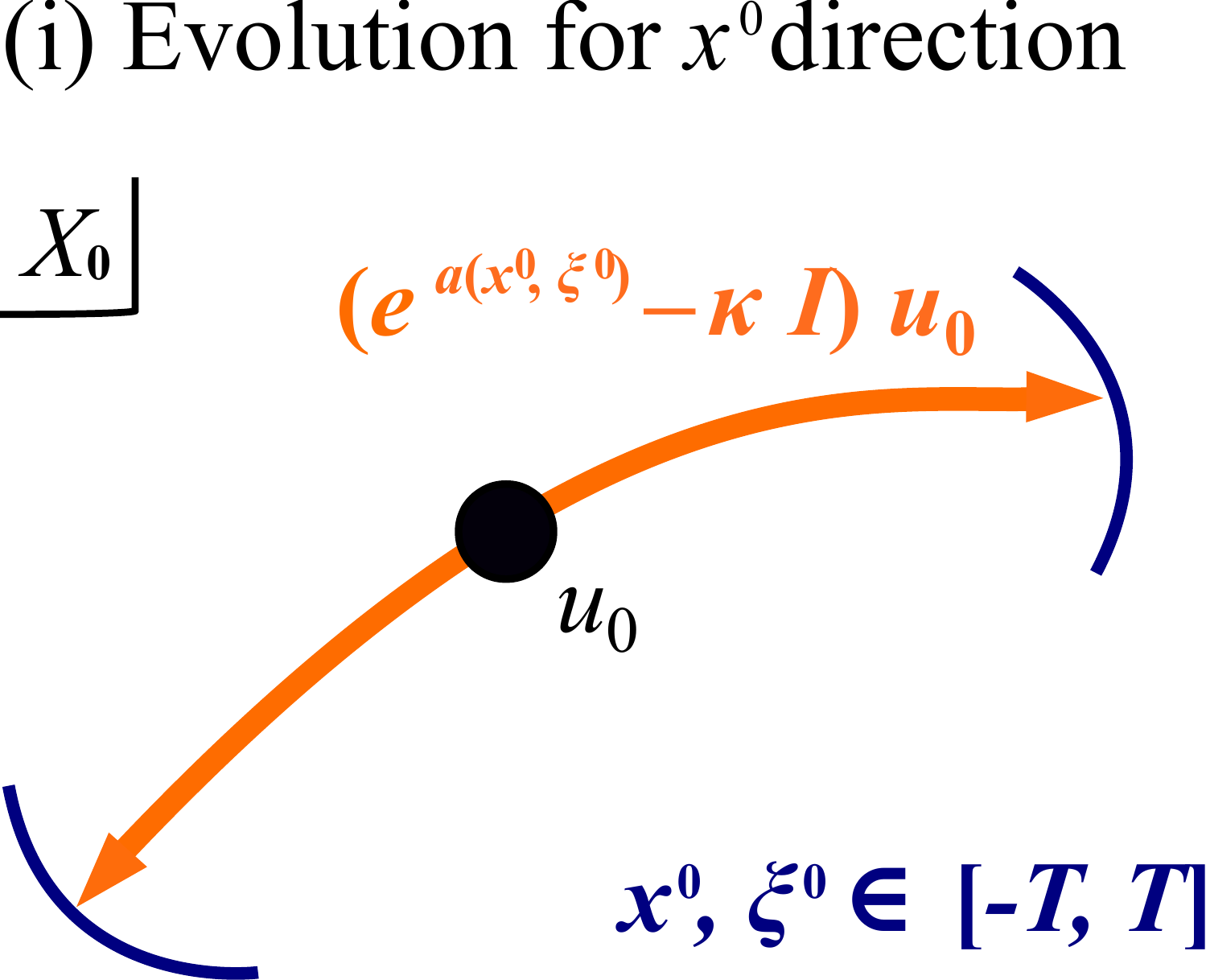} \hspace{6mm}
\includegraphics[width=48mm]{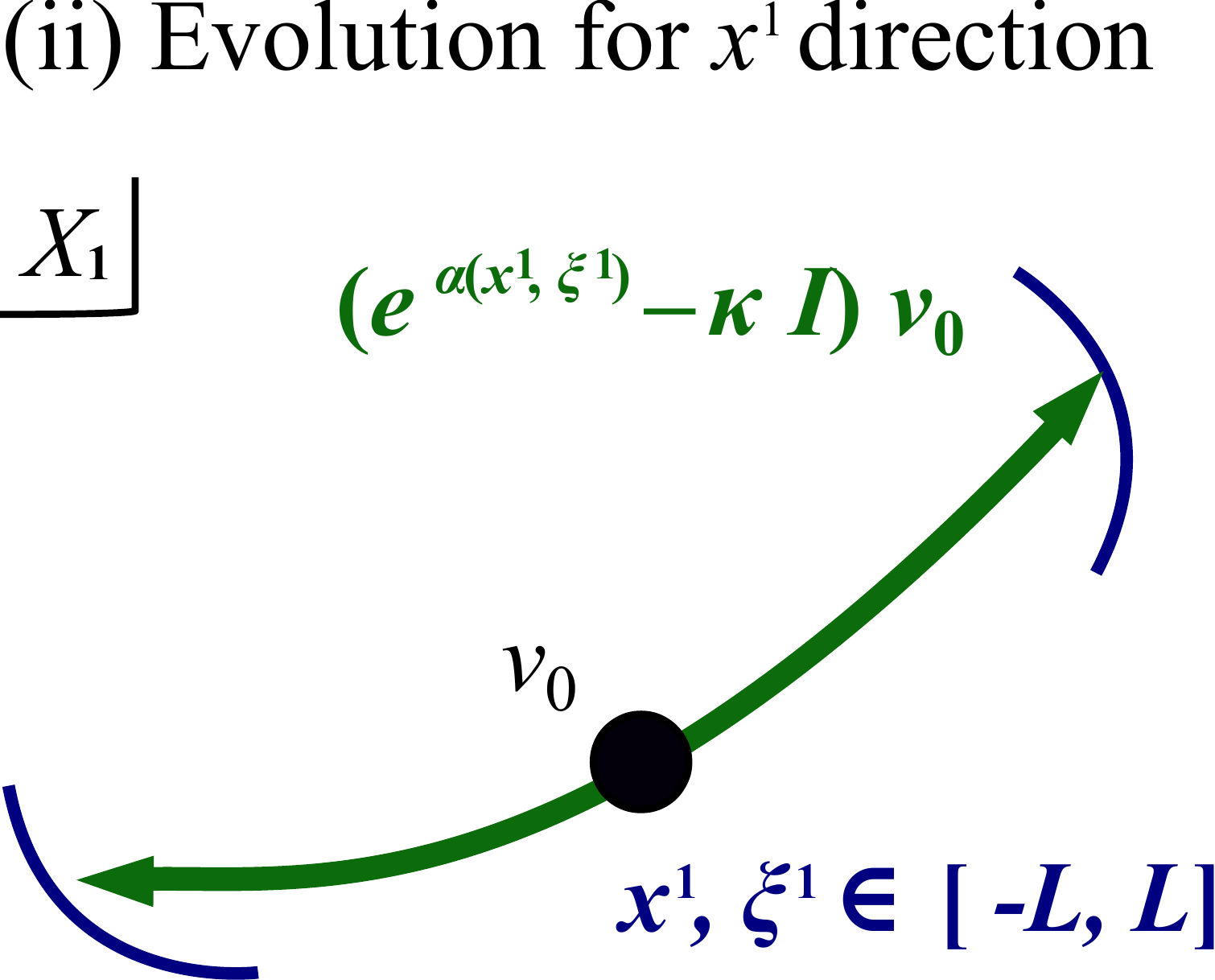} \hspace{6mm} 
\includegraphics[width=48mm]{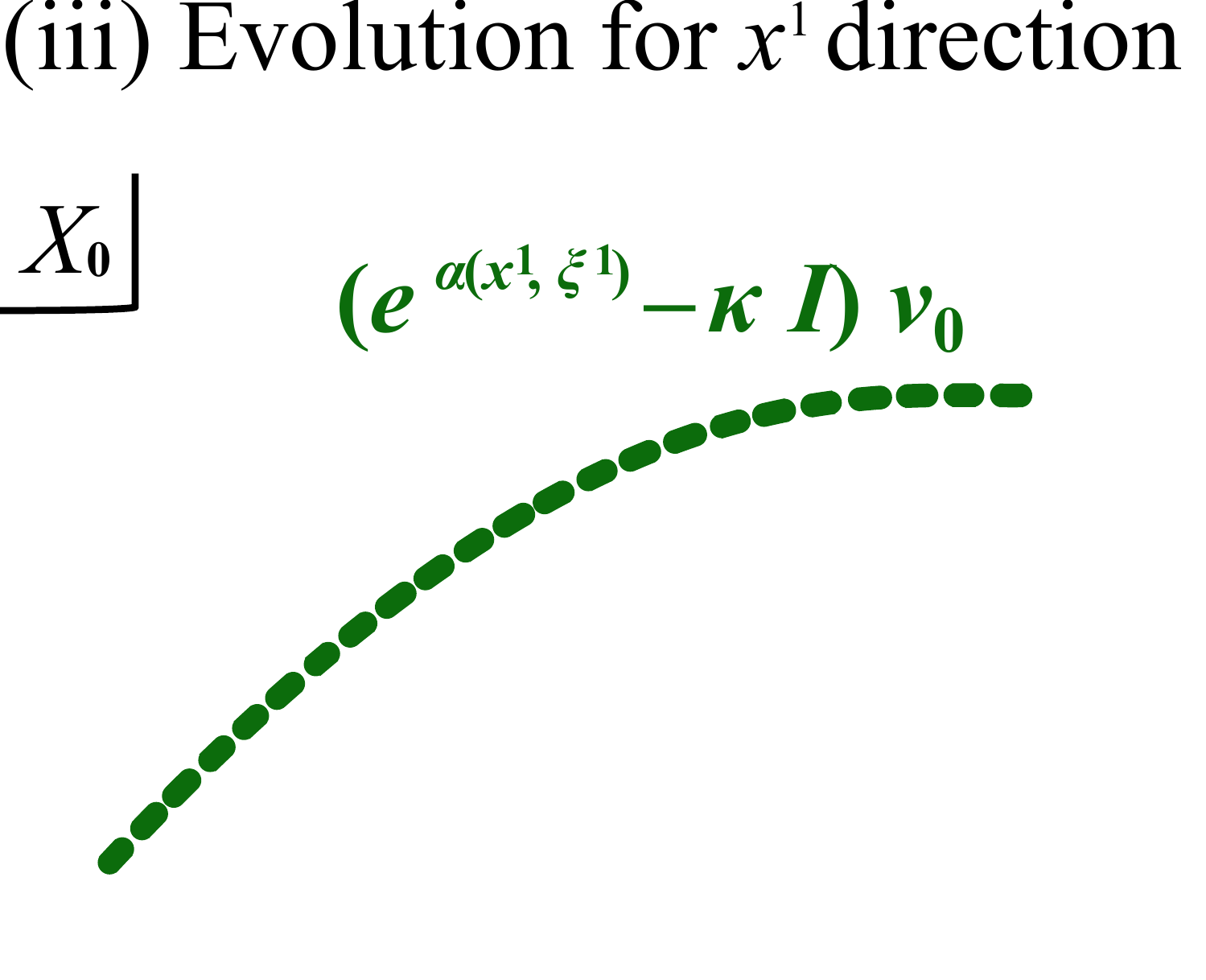} \vspace{6mm} \\
\caption{ \label{fig2} 
(Color online) Three stages for obtaining discrete trajectory by means of the relativistic formulation.
In panel (i)  $e^{a(x^0,\xi^0)} - \kappa I$ acting on $u_0 \in X_0$ is shown as a trajectory in $X_0$. 
In panel (ii) $e^{\alpha(x^1,\xi^1)} - \kappa I$ acting on elements in $X_1$ is obtained by changing the evolution direction. 
In panel (iii)  the regularized trajectory $e^{\alpha(x^1,\xi^1)} - \kappa I v_0 \in C^0(-L,L;L^2(-T,T))$ is regarded as an element in $L^2(-T,T;C^0(-L,L)) $ that leads to the concept of generating discrete groups.
}
\end{figure}

In conclusion, the logarithmic representation is obtained for the infinitesimal generators of group and semigroups of operators.
By associating the logarithmic representation with a relativistic formulation, some non-trivial properties are noticed.
In particular the following properties are remarkable for the logarithmic representation of relativistic evolution equations.  \vspace{3mm} \\
%%%%%%
\underline{\bf Breaking and recovery of time reversal symmetry} ~
%\subsection{Logarithmic representation of infinitesimal generators for general evolution families}
The representation similar to (\ref{logex}) has been originally obtained for the invertible evolution families~\cite{17iwata-1}, while it is generalized to non-invertible evolution families in this article.
%%%%%%%
Indeed, as shown in the above proof, the invertible property of $U(x^i,\xi^i)$ is not used.
Under the validity of boundedness of $U(x^i,\xi^i)$ on $X$, the removal of invertible criterion is essentially realized by the introduction of nonzero $\kappa \in {\mathbb C}$ with a sufficiently large amplitude.
On the other hand, the indispensable conditions for obtaining this kind of logarithmic representations are the boundedness of the spectral set of $U(x^i,\xi^i)$ and the commutation assumption, where the bounded interval $-L \le x^i,\xi^i \le L$ is also necessary.
%The boundedness is satisfied by the infinitesimal generators of semigroups and groups defined in the standard theory of abstract evolution equations, and the commutation is satisfied at least by $x^i$-independent infinitesimal generators. 
As a result, the concept of regularized trajectory being related to analytic semigroups (for a textbook, see Ref.~\cite{79tanabe}) follows (for the relation, see Ref.~\cite{17iwata-1}).
%%%%%%
The regularized trajectory is a local-in-$x^i$ representation, so that the locality is also essential to the representation.
Furthermore the regularization is associated with the recovery of local time reversal symmetry if $x^i$ is equal to $x^0$. 
It should bring about one mathematical clue to understand the breaking mechanism of time reversal symmetry. 
\vspace{3mm} \\
%%%%%%
\underline{\bf Nonlinearity} ~
The operator $U(x^i,\xi^i)$ can be either linear or nonlinear semigroup.
Indeed, the linearity of the semigoup is not used in the proof.
As shown in Ref.~\cite{18iwata}, the nonlinearity of semigroup can appear simply by altering the evolution direction under a suitable identification between the infinitesimal generator and the evolution operator.
In this sense, Eq.~(\ref{k-eq1}) is regarded as a local-in-$x^i$ linearized equation, if $U(x^i,\xi^i)$ is a nonlinear semigroup (semigroup related to the nonlinear equations).
That is, the logarithmic representation provides a way to obtain a linearized infinitesimal generator from nonlinear semigroups. 
%%%%%%
As suggested in Ref.~\cite{18iwata}, the relation between evolution operator and its infinitesimal generator is essentially similar to the Cole-Hopf transform, and the similarity is seen by introducing the relativistic formulation.
It should bring about a key to understand the emergence of a certain kind of nonlinearity. 
The present result can be used to linear evolution equations of hyperbolic type and the related quasi-linear evolution equations~\cite{75kato}.
\vspace{3mm} \\
%%%%%%
\underline{\bf Discrete property} ~
For example, in case of two-dimensional space-time distribution, let the $C^0$-semigroup for $x^0$ direction exist for a Cauchy problem:
\[ \begin{array}{ll}
\partial_{x^0} U(x^0, \xi^0)~u_0  = K(x^0) U(x^0,\xi^0) ~u_0,  \vspace{2.5mm} \\
 \partial_{x^0} a(x^0,\xi^0) = \partial_{x^0} {\rm Log}  (U(x^0,\xi^0) + \kappa I),
\end{array} \]
in $X_0 := L^2(-L,L)$, where a certain complex number $\kappa$ is taken from the resolvent set of $U(x^i,\xi^i)$.
The corresponding dynamical system is illustrated in panel (i) of Fig.~\ref{fig2}.
Furthermore let the same equation be written as
\[ \begin{array}{ll}
\partial_{x^1} V(x^1, \xi^1)~v_0  = {\mathcal K}(x^1) V(x^1,\xi^1) ~v_0,  \vspace{2.5mm} \\
 \partial_{x^1} \alpha(x^1,\xi^1) = \partial_{x^1} {\rm Log}  (V(x^1,\xi^1) + \kappa I),
\end{array} \]
in $X_1:= L^2(-T,T)$, and the corresponding dynamical system is illustrated in panel (ii) of Fig.~\ref{fig2}.
In this situation, using $\alpha(x^1,\xi^1)$ instead of $a(x^0,\xi^0)$, the corresponding dynamical system contains a discrete trajectory in $X_0$ (the right panel of Fig~\ref{fig2}).
Indeed, the trajectory is $L^2$-function with respect to $x^0$, and $C^0$ function with respect to $x^1$. 
That is, the relativistic treatment leads to the discrete evolution (for a theory including the discrete evolution, see the variational method of abstract evolution equation \cite{61lions}).
The discrete evolution to the $t$-direction ($x^0$-direction) is useful to analyze the stochastic differential equations within the semigroup theory of operators.

% Acknowledgement
\section{ACKNOWLEDGMENTS}
The autuhor is grateful to Prof. Emeritus Hiroki Tanabe of Osaka University.
This work was supported by JSPS KAKENHI Grant No. 17K05440.
% References

\nocite{*}
\bibliographystyle{aipnum-cp}%

\bibliography{sample}%

\end{document}